\documentclass[12pt,reqno]{article}

\usepackage[usenames]{color}
\usepackage{amssymb}
\usepackage{amsmath}
\usepackage{amsthm}
\usepackage{amsfonts}
\usepackage{amscd}
\usepackage{graphicx}
\usepackage{xcolor}
\usepackage{float}

\usepackage[colorlinks=true,
linkcolor=webgreen,
filecolor=webbrown,
citecolor=webgreen]{hyperref}

\definecolor{webgreen}{rgb}{0,.5,0}
\definecolor{webbrown}{rgb}{.6,0,0}

\usepackage{fullpage}

\usepackage{graphics}
\usepackage{latexsym}
\usepackage{epsf}
\usepackage{breakurl}

\newcommand{\seqnum}[1]{\href{https://oeis.org/#1}{\rm \underline{#1}}}

\theoremstyle{plain}
\newtheorem{theorem}{Theorem}

\newtheorem{lemma}[theorem]{Lemma}
\newtheorem{proposition}[theorem]{Proposition}

\theoremstyle{definition}

\newtheorem{conjecture}[theorem]{Conjecture}
\newtheorem{question}[theorem]{Question}

\theoremstyle{remark}

\DeclareMathOperator{\rad}{rad}

\begin{document}

\begin{center}
\vskip 1cm{\LARGE\bf Must a Primitive Non-Deficient Number Always Have a Small Component?}
\vskip 1cm
\large
Joshua Zelinsky\\
Department of Mathematics\\ 
Hopkins School\\ 
New Haven, CT 06515\\
USA\\
\href{mailto:jzelinsky@hopkins.edu}{\tt jzelinsky@hopkins.edu} \\
\end{center}

\vskip .2 in
\begin{abstract} Let $n$ be a primitive non-deficient number where $n=p_1^{a_1}p_2^{a_2} \cdots p_k^{a_k}$ where $p_1, p_2, \ldots, p_k$ are distinct primes. We prove that there exists an $i$ such that $$p_i^{a_i+1} < 2k(p_1p_2p_3\cdots p_k).$$ We conjecture that in fact one can always find an $i$ such that $p_i^{a_i+1} < 2p_1p_2p_3\cdots p_k$.
\end{abstract}

\section{Introduction}

Let $\sigma_d(n)$ be the sum of the $d$th powers of divisors of $n$. We  write $\sigma(n)$ for $\sigma_1(n)$.
Given a positive integer $n$, we write $h(n) = \sigma_{-1}(n) = \frac{\sigma(n)}{n}$. We write 
 $$H(n) =\prod_{p \mid n, \, p \mathrm{\,\, prime}} \frac{p}{p-1}.$$ In what follows, we write $n= p_1^{a_1}p_2^{a_2} \cdots p_k^{a_k}$ where $p_1, p_2, \ldots, p_k$ are distinct primes and the $a_i$ are positive integers. We write $\rad (n)$ to be the radical of $n$. That is, $\rad (n) = p_1p_2 \cdots p_k$. We also set $R=\rad (n)$.

Recall, a number $n$ is said to be deficient if $\sigma(n) < 2n$. Deficient numbers form \seqnum{A005100}. If $n$ is non-deficient, then so is any multiple of $n$.  Thus, there is an interest in non-deficient numbers which are primitive in the sense that such a  number is non-deficient and also has no non-deficient divisors.  Primitive non-deficient numbers are sequence \seqnum{A006039} in the OEIS. They have been studied since at least Dickson \cite{Dickson} who proved that for any fixed $k$ there are only finitely many odd primitive non-deficient numbers. An asymptotic for the number of primitive non-deficient numbers which are at most $x$ is still open. The best known bounds are due to Avidon \cite{Avidon}. In a different direction, the non-deficient numbers have a positive density and there are estimates of their natural density.  Kobayashi \cite{Kobayashi}, showed that their natural density is at least $0.24760444$, while Del\'eglise  \cite{Deleglise} obtained an upper bound of $0.248$. More recently, a preprint by McNew and Setty \cite{McNewSetty} obtained that their natural density is at least $0.247619608$ and at most $0.247619658$.

It is a classical observation dating to the 19th century that the functions $h(n)$ and $H(n)$ are closely linked. In particular, $h(n) \leq H(n)$ with equality if and only if $n=1$. Moreover, $H(n)$ is the best upper bound one can have on $h(n)$ that depends only on the distinct prime factors of $n$. In particular,

$$\lim_{t \rightarrow \infty} h(n^t) = H(n) = H(\rad(n)).$$

We tentatively suspect that

\begin{conjecture}\label{mainconjecture} If $n$ is a primitive non-deficient number with $n={p_1}^{a_1}{p_2}^{a_2} \cdots {p_k}^{a_k}$, 
and $p_1 < p_2 < \cdots < p_k$,
then there exists an $i$, $1 \leq i \leq k$ such that $p_i^{a_i+1} < 2R$.
\end{conjecture}

Conjecture \ref{mainconjecture} is false if $2R$ is replaced by $R$, although the author knows of only one counterexample, namely, $n=4164647056479 =  3^7 7^4  13^3  19^2 $, where  $R=5187$, and the smallest of $p_i^{a_i+1}$ is $3^8=6561$. It may be that this is the only counterexample. A slightly weaker version of Conjecture \ref{mainconjecture} would replace $2R$ with $p_1R$, and we are more confident in that conjecture. 

We prove the following weaker result:

\begin{theorem} Let $n$ be a primitive non-deficient number with $n={p_1}^{a_1}{p_2}^{a_2} \cdots {p_k}^{a_k}$ with $p_1 < p_2 < \cdots < p_k$. Let $R=\rad(n) = p_1p_2\cdots p_k$. Then there is an $i$ such that \label{maintheoremboundwithkR}
\begin{equation}\label{Generalinequalityforprimitivenondeficientnumbers}p_i^{a_i+1} < 2kR.\end{equation}
\end{theorem}

It is plausible that the right-hand side of Inequality \eqref{Generalinequalityforprimitivenondeficientnumbers} can be replaced with $p_1R$, which would still be weaker than Conjecture \ref{mainconjecture}.

Conjecture \ref{mainconjecture} and Theorem \ref{maintheoremboundwithkR} are motivated by three interrelated lines of thought. 

First, we have the following intuition: if $n$ is a non-deficient number (primitive or not), then one must have $2 \leq h(n) < H(n)$, but if all of the $a_i$ are ``large'' then for each $i$, $h(p_i^{a_i})$ will be close to $H(p_i)$, and thus $h(n)$ will be very close to $H(n)$. Thus, if the $a_i$ are large enough, then for at least one $p_j$, $h(\frac{n}{p_j})$ will still be very close to $H(n)$, and thus $\frac{n}{p_j}$ will be non-deficient (since it will be very close to $H(n)$) and so $n$ will not be primitive. Thus, if the $p_i$ are all fixed and $n$ is in fact primitive, then at least one of the $p_i^{a_i}$ must be small. Conjecture \ref{mainconjecture} and Theorem \ref{maintheoremboundwithkR} are essentially making this intuition precise. The radical of $n$ is a natural function to make this notion precise since the radical exactly contains the information of what the $p_i$ are.  

Second, one concern in some recent work concerning odd perfect numbers has been the graph representation of the components of the odd perfect number. In particular, given an odd perfect number $n =p_1^{a_1}p_2^{a_2} \cdots p_k^{a_k}$, we can construct a weighted directed graph $G$ with $k$ vertices labeled $p_1, p_2, \ldots, p_k$ with an edge from $p_i$ to $p_j$ of weight $m$ where $p_i^{m} \mid \mid \sigma(p_j^{a_j})$. Recent papers have used this graph framework or similar frameworks.  One example is a paper Bibby, Vyncke, and the author \cite{BVZ}  which proved new upper bounds on the size of the third largest prime factor of an odd perfect number. A similar graph approach was used in a paper of Yamada \cite{Yamada2}  to produce new results about solutions to the equation 

\begin{equation} \sigma(n) = \rad(n)^2 \label{DeKoninck}, \end{equation}
where it is currently unknown if the only solutions are $n=1$ and $n=1782$. Other notable work on solutions to Equation \eqref{DeKoninck} also has been done by Broughan, De Koninck, K\'{a}tai, and Luca, \cite{BKKL} with related work by Luca \cite{Luca}. In this context, if $n$ is an odd perfect number, then Conjecture \ref{mainconjecture} and Theorem \ref{maintheoremboundwithkR} can be thought of as bounds on the out-degree of vertices in this graph. 

The third line of thought also concerns odd perfect numbers. Over the last few years, multiple papers have developed inequalities relating an odd perfect number $n$ to $\rad(n)$. Acquaah and Konyagin \cite{AcquaahKonyagin} proved that $R < 2n^\frac{2}{3}$, and Luca and Pomerance \cite{LucaPomerance} proved that $R < 2n^\frac{17}{26}$. Subsequently, Klurman \cite{Klurman} proved that there is a constant $c$ (independent of $n$) such that $R < cn^\frac{9}{14}$, and it is not hard to show that one may take $c=2$ in Klurman's argument. Ellia \cite{Ellia} showed that under some reasonable technical assumptions one must have $R < n^\frac{1}{2}$, and Ochem and Rao \cite{OchemRaoradica} extended Ellia's results further. 

Now, all of these results are upper bounds for $R$ in terms of $n$. But one can also prove upper bounds for $n$ in terms of $R$. In particular, Nielsen \cite{Nielsen} proved that if $n$ is an odd perfect number then 
\begin{equation}
    n < R^{2^k-1} \label{Nielsenupperbound}.
\end{equation}

An immediate consequence of Equation \eqref{Nielsenupperbound} is that if $n$ is an odd perfect number then there exists an $1 \leq i \leq k$ such that $p_i^{a_i} \mid \mid n$, and where
\begin{equation}
    p_i^{a_i+1} < R^{\frac{2^k}{k}}. \label{Nielsenboundcorrolary}
\end{equation}
Thus, Conjecture \ref{mainconjecture} and Theorem \ref{maintheoremboundwithkR} can be thought of as generalizations of Equation \eqref{Nielsenboundcorrolary}, in that the inequality is tighter and applies also to all primitive non-deficient numbers. 

Now, with all of these motivations, a reader may reasonably wonder why look at an inequality involving $p_i^{a_i+1}$ and $R$ rather than $p_i^{a_i}$ and $R$? There are two answers for this, one somewhat unsatisfying and one hopefully more satisfying. The less satisfying answer is that the techniques used in this paper naturally seem to give inequalities involving $p_i^{a_i+1}$. The more satisfying answer is that these inequalities are equivalent to an inequality for $p_i^{a_i}$ in terms of all the other prime factors of $n$. That is, if one has an inequality of the form

$$p_i^{a_i +1} < f(k)R,$$ for some function $f(k)$ then this is the same as 

$$p_i^{a_i} < f(k)\prod_{j=1, j\neq i}^{k} p_i,$$

so this translates to an inequality for a component in terms of the product of all the other prime factors.

\section{Main results}
We first recall the following basic fact about non-deficient numbers which has a straightforward proof. 
\begin{lemma} Let $n$ be a number of the form $n=2^ap$ where $p$ is an odd prime. Then $n$ is a primitive non-deficient number if and only if $2^a <p \leq 2^{a+1} -1$. Also, all odd numbers with exactly two distinct prime divisors are deficient.
\label{Classificationofprimitivenondeficientwhenkis2}
\end{lemma}
One consequence of Lemma \ref{Classificationofprimitivenondeficientwhenkis2} we prove is that Theorem \ref{maintheoremboundwithkR} holds whenever $n$ is even and $i=1$. 

\begin{lemma}   Let $n$ be an odd primitive non-deficient number with $n={p_1}^{a_1}{p_2}^{a_2} \cdots {p_k}^{a_k}$, 
and $p_1 < p_2 < \cdots < p_k$. Assume further that $k=3$ or $k=4$
then there exists an $i$, $1 \leq i \leq k$ such that $p_i^{a_i+1} < 2R$.
 \label{Checkforexactly3or4distinctprimefactors}
\end{lemma}
\begin{proof} There are no odd non-deficient numbers with fewer than 2 distinct prime factors, and the only odd primitive non-deficient numbers with 3 distinct prime factors are with a little computation seen to be just the elements of the set  $$\{945, 1575, 2205, 7425, 78975, 131625, 342225, 570375\}$$ and one can check that all elements of this set satisfy the desired inequality. The set for odd primitive non-deficient numbers with 4 distinct prime factors is longer, but this is still a straightforward computation.
\end{proof}

The next lemma is old, with a slightly weaker version due to Servais \cite{Servais}. Although this lemma is often phrased just for odd perfect numbers, the lemma applies to any non-deficient number. 

\begin{lemma} Let $n$ be a non-deficient number with $n= p_1^{a_1}p_2^{a_2} \cdots p_k^{a_k}$ with primes $p_1 < p_2 < \cdots < p_k$. Then

\begin{equation} \frac{3}{2}p_1 -2 <  k. \label{Grunequation}
\end{equation}

Additionally, if $p_1 \geq 5$, then 

\begin{equation} p_1 +2 \leq k. \label{Grunlikeneededotherinequality}
\end{equation}

\label{Grun}
\end{lemma}
\begin{proof} Inequality \eqref{Grunequation} is essentially just the argument from Gr{\"u}n  \cite{Grun}.

The second inequality follows immediately from   Inequality \eqref{Grunequation} for $p_1 \geq 11$, so it only needs to be verified for $p_1=5$ and $p_1=7$.

We verify it for $p_1=5$. The proof for $p_1=7$ is nearly identical. Assume $n$ is a non-deficient number with $n= p_1^{a_1}p_2^{a_2} \cdots p_k^{a_k}$ with primes $p_1 < p_2 < \cdots < p_k$, and that $p_1=5$. Assume further that $k \leq 6$. Since $n$ is non-deficient, we have

\begin{equation}
    2 \leq h(n) < H(n) = \prod_{i=1}^k \frac{p_i}{p_i-1}. \label{standardhlessthanHbound}
\end{equation}
Now, since $\frac{x}{x-1}$ is a decreasing function and $k \leq 6$, and $p_1=5$, we have from Equation \eqref{standardhlessthanHbound},

\begin{equation}
    2 < \frac{5}{4}\frac{7}{6}\frac{11}{10}\frac{13}{12}\frac{17}{16}\frac{19}{18} < 2,
\end{equation}

which is a contradiction. So $k \geq 7$.
\end{proof}

Note that a stronger version of Inequality \eqref{Grunequation} can be found in a prior paper by the author \cite{Zelinskybig}, and similarly tighter bounds have been found by others \cite{AslaskenKirfel, Stone}  but we do not need it here. We have then the following result for even primitive non-deficient numbers. 

\begin{proposition} Let $n$ be an even primitive non-deficient number, with $R=\rad(n)$. Let $a$ be the integer such that $2^a \mid \mid n$. Then $2^{a+1} < R$. Thus, every even primitive non-deficient number satisfies the radical inequality. \label{evenprimitivenondeficientshaveatightbound}
\end{proposition}
\begin{proof}  Let $n$ be an even primitive non-deficient number. Let $a$ be the positive integer such that $2^a \mid \mid n$. If $n$ were a power of 2, it would be deficient, so $n$ must have at least one odd prime divisor. Call that divisor $p$. Consider then $m=2^a p$. Since $n$ is a primitive non-deficient number, it must be the case that $m$ is deficient. Then, from Lemma \ref{Classificationofprimitivenondeficientwhenkis2}, we must have $2^a <p$, so $2^{a+1} < 2p \leq R$.
\end{proof}

Even primitive non-deficient numbers are much more common than odd primitive non-deficient numbers, so we can say that Conjecture \ref{mainconjecture} is at least true for most primitive non-deficient numbers. We prove a weaker version of the conjecture but we first need a few lemmas.

\begin{lemma} If $k$ and $m$ are real numbers such that $k \geq 2$ and $0 < m < \frac{1}{k(k+1)}$, then 

\begin{equation}\frac{1}{1-km} < 1 + (k+1)m. \label{KMinequality}\end{equation} \label{KMlemma}
\end{lemma}
\begin{proof} Assume that $k$ and $m$ are real numbers such that $k \geq 2$ and $0 < m < \frac{1}{k(k+1)}$. Then, since $km<1$, we have that $1-km$ is positive, and so our desired inequality is equivalent to showing that 
$$1 < (1-km)(1+(k+1)m) = 1 + m(1 - k(k + 1)m), $$ and thus our inequality is equivalent to the claim that $m(1 - k(k + 1)m)$ is positive, but this is true since $m>0$ and $1 - k(k + 1)m>0$, since $m < \frac{1}{k(k+1)}$.
\end{proof}

\begin{lemma} Let $n$ be a positive integer, and let $p$ be a prime such that $p^a \mid \mid n$. Then $$\frac{h(n)}{h(\frac{n}{p})} \leq 1 + \frac{1}{p^a}. $$
\label{correctedPAlemma}
\end{lemma}

\begin{proof} We break this into two cases, depending on whether $a=1$ or $a \geq 2$. If $a=1$, since $h(n)$ is multiplicative we have, 
$$\frac{h(n)}{h(\frac{n}{p})}  = h(p)=1+ \frac{1}{p}.$$
Now, assume that $a \geq 2$. Then we have

$$\frac{h(n)}{h(\frac{n}{p})} = \frac{1 + \frac{1}{p} \cdots + \frac{1}{p^{a-1}}+\frac{1}{p^a}}{1 + \frac{1}{p} \cdots \frac{1}{p^{a-1}}} = 1 + \frac{\frac{1}{p^a}}{1 + \frac{1}{p} \cdots \frac{1}{p^{a-1}}} < 1 + \frac{1}{p^a}.$$
\end{proof}

\begin{lemma}\label{GeneralizedPuchtalemma} Let $n$ be a primitive non-deficient number with $n=p_1^{a_1}p_2^{a_2} \cdots p_k^{a_k}$ with $p_1 < p_2 < \cdots < p_k$. Assume further that $(p_1,k) \neq (3,3)$. Let $R=\rad(n) = p_1p_2\cdots p_k$. Assume further that $$H(n) \geq 2 + \alpha$$ for some positive $\alpha <1$. Then there exists an $i$, $1 \leq i \leq k$ such that 
$$p_i^{a_i+1} < \max\left(\frac{2(k+2+p_1)}{\alpha },k(k+1)\right).$$
\end{lemma}
\begin{proof} This proof is essentially an adaptation of part of the argument used in prior work by Heath-Brown \cite{Heath-Brown}. Assume $n$ is a primitive non-deficient number with $n=p_1^{a_1}p_2^{a_2} \cdots p_k^{a_k}$ with $p_1 < p_2 < \cdots < p_k$ and let $R=\rad(n) = p_1p_2\cdots p_k$. Assume  that $$H(n) \geq 2 + \alpha$$ for some positive $\alpha <1$. So we have $$H(n) = \prod_{i=1}^k \frac{p_i}{p_i -1}   \geq 2 + \alpha .$$  We also have that
 $$h(n) = H(n) \prod_{i=1}^k \left(1 - \frac{1}{p_i^{a_i+1}}\right) = \prod_{i=1}^k \left(\frac{p_i}{p_i-1}\right)\left(1 - \frac{1}{p_i^{a_i+1}}\right).$$

But we also have $$\prod_{i=1}^k \left(\frac{p_i}{p_i-1}\right)\left(1 - \frac{1}{p_i^{a_i+1}}\right) \geq \left(\prod_{i=1}^k \frac{p_i}{p_i-1}\right)\left(1 - \sum \frac{1}{p_i^{a_i+1}}\right),$$
and $$ \left(\prod_{i=1}^k \frac{p_i}{p_i-1}\right)\left(1 - \sum \frac{1}{p_i^{a_i+1}}\right) \geq \left(2+ \alpha \right)\left(1 - k \max  \frac{1}{p_i^{a_i+1}}\right).$$

Now, reducing any of the $a_i$ by 1 or more to get a new number $n'$ yields a deficient number since $n$ is primitive non-deficient. Thus for any prime $p_i$ we have $$h(n/p_i) < 2,$$ which we combine with   Lemma \ref{correctedPAlemma}  to conclude that for any $p_i$, 
$$h(n) < 2\left(1 + \frac{1}{p_i^{a_i}}\right). $$
We have then,

$$2\left(1 + \min \frac{1}{p_i^{a_i}} \right)> (2+\alpha )\left(1 - k \max  \frac{1}{p_i^{a_i+1}}\right). $$
Set $m= \max \frac{1}{p_i^{a_i+1}}$. Note that $m = \frac{1}{\min p_i^{a_i +1}}$ and that 
$\min \frac{1}{p_i^{a_i}} \leq p_1 m$. Thus, 
$$\frac{1+p_1m}{1-km} > 1+\frac{\alpha}{2}. $$
If one has a concern about sign issues in the above, note that if $km  \geq 1$, then one has that $\min p_i^{a_i+1} < k$, which is stronger than our desired inequality, so we can assume that $km  < 1$. Note that $k \geq p_1$ by Lemma \ref{Grun}, and the fact that no non-deficient number is a power of 2. 

Now, we split into two cases, where $m < \frac{1}{k(k+1)}$ and where $m \geq \frac{1}{k(k+1)}$. Notice that since $k(k+1)$ cannot be a power of a prime, these two quantities cannot be equal.

Consider the case where $m > \frac{1}{k(k+1)}$. This means that we have $$\min p_i^{a_i+1} < k(k+1),$$ which is one of the terms in our maximum needed.  If $$m < \frac{1}{k(k+1)},$$ then by Lemma \ref{KMlemma}, we have 

$$\frac{1}{1-km} < 1 + (k+1)m.$$
Hence,

\begin{equation}\frac{1+p_1m}{1-km} < (1+ mp_1)( 1 + (k+1)m) = 1 + (k+1)m + mp_1 + m^2p_1(k+1). \label{preludeinequality} \end{equation}

Now, we use that $p_1 < k$  with Inequality \eqref{preludeinequality} to obtain the following:

$$1 + (k+1)m + mp_1 + m^2p_1(k+1) < 1 + m(k+1+p_1) + m^2k^2 < 1 + m(k + 2 +p_1),$$

and so $1 + m(k + 2 +p_1) > 1+ \frac{\alpha}{2} $,
which yields$$ m(k + 2 +p_1) > \frac{\alpha}{2}. $$
Thus, 
 \begin{equation}\label{copyofstrongermininequalityinlemmaform} \min p_i^{a_i+1} < 2\frac{(k+2+p_1)}{\alpha }.\end{equation}
 
\end{proof}

We are now ready to prove Theorem \ref{maintheoremboundwithkR}.
\begin{proof}
 We may, based on the earlier results, assume that $p_1 \geq 3$. By Lemma \ref{Checkforexactly3or4distinctprimefactors} we may also assume that $k \geq 5$.
 
 Since $n$ is non-deficient, we have 
$$2 \leq h(n) < H(n).$$
We have $$H(n) = \prod_{i=1}^k \frac{p_i}{p_i -1}  \geq 2 + \frac{1}{(p_1-1)(p_2-1) \cdots (p_k-1)} \geq 2 + \frac{2}{\prod_{i=1}^k p_i-1} > 2 + \frac{2}{R}.$$  
Note that on the last step in the above chain of inequalities, we are using that $\prod_{i=1}^k \frac{p_i}{p_i-1} >2 $ so $(p_1-1)(p_2-1) \cdots (p_k-1) < \frac{R}{2}$.

We may then apply Lemma  \ref{GeneralizedPuchtalemma} with $\alpha = \frac{2}{R}$. Thus, there exists an $i$ such that  \begin{equation}\label{copyofstrongermininequalitywithmaxterm} p_i^{a_i+1} < \max (R(k+2+p_1), k(k+1)).\end{equation}
The first term in the above maximum is clearly larger than the second term  so we have 

\begin{equation}\label{copyofstrongermininequalitywithoutmaxterm} p_i^{a_i+1} < R(k+2+p_1).\end{equation}
We now apply Lemma \ref{Grun} along with the fact that $k \geq 5$ to get
$p_1 +2 \leq k$ 

and thus $$R(k+2+p_1) \leq 2Rk. $$

We thus have 

$$p_i^{a_i+1} \leq 2Rk.$$ But  since the left-hand side is odd and the right-hand side is even, the inequality must be strict and so we are done.
\end{proof}

 Note that one can get a slightly tighter, but substantially less nice looking bound by replacing Lemma \ref{Grun} in the above argument with the sort of inequality used in section 5 of the author's aforementioned paper  \cite{Zelinskybig}.

\section{Related results and further conjectures}
How tight is Conjecture \ref{mainconjecture}?  Some of our odd primitive non-deficient numbers with 3 prime factors just barely manage to satisfy this inequality. Some of worst case situations occur when $n = 2^ap$ and $p$ is a Fermat prime. That is, we have $p=2^a +1$. In this scenario, $2^{a+1} = R-2$ and so these are close to our inequality being violated. However, the standard conjecture is that there are only finitely many Fermat primes. 

The main theorem of this paper and Conjecture \ref{mainconjecture} motivate the following definitions and additional conjectures. 

Let $A$ be the set of numbers satisfying $p_i^{a_i+1}< 2R$ for some $i$. It is not hard to show that such numbers have density 1. Let $B$ be the set of numbers satisfying $p_i^{a_i+1} < 2R$ for all $i$. An example of such a number is $105$. Note that there are some primitive non-deficient  numbers in the set $B$. For example,  $945$ is in $B$. What can we say about density of numbers from $B$? Note that $B$ has natural density zero due to classical results about the distribution of distinct prime factors such as those by Billingsley \cite{Billingsley}. Given a set $S$, we write $S(x)=|\{n \in S, n \leq x\}|$. It seems natural to ask for an asymptotic for $B(x)$. 

Although 945 is in $B$, such examples seem to be rare.  Let $C$ be the set of primitive non-deficient numbers.  Let $C_O$ be the set of odd primitive non-deficient numbers. Then we suspect that
\begin{conjecture}
$$\lim_{x \rightarrow \infty} \frac{(B \cap C)(x)}{B(x)}=0.$$
\end{conjecture}
Similarly, it seems likely that
\begin{conjecture}
$$\lim_{x \rightarrow \infty} \frac{(B \cap C)(x)}{C(x)}=0,$$
\end{conjecture}

We have similar conjectures for $C_O$. In particular, we have the following:

\begin{conjecture}
$$\lim_{x \rightarrow \infty} \frac{(B \cap C_O)(x)}{B(x)}=0.$$
\end{conjecture}
\begin{conjecture}
$$\lim_{x \rightarrow \infty} \frac{(B \cap C_O)(x)}{C_O(x)}=0.$$
\end{conjecture}

We have the following generalization: Let $M$ be a set of primes, and let $C_M$ be the set of primitive non-deficient numbers which do not have any element of $M$ as a prime divisor. So, $C_O = C_{\{2\}}$. Then the following seem plausible.

\begin{conjecture} Let $M$ be a non-empty finite set of primes. Then
$$\lim_{x \rightarrow \infty} \frac{(B \cap C_M)(x)}{B(x)}=0.$$ \label{StrongversionforM1}
\end{conjecture}
\begin{conjecture} Let $M$ be a non-empty finite set of primes. Then \label{StrongversionforM2}
$$\lim_{x \rightarrow \infty} \frac{(B \cap C_M)(x)}{C_M(x)}=0.$$
\end{conjecture}

We can also apply Lemma \ref{GeneralizedPuchtalemma} to obtain the following result, which essentially says that if $n$ is an odd primitive non-deficient number, then either its largest prime factor is small compared to the radical or it has a component which is small compared to the radical.

\begin{proposition}  Let $n$ be an odd  primitive non-deficient number with $n=p_1^{a_1}p_2^{a_2} \cdots p_k^{a_k}$ with $p_1 < p_2 < \cdots < p_k$.
Either \begin{equation}p_k < \sqrt{R}\label{pklessthansquarerootofR}\end{equation} or there exists an $i$ such that \begin{equation}p_i^{a_i+1} < (k+2+p_1)(p_k-1).\label{mincomponentnotmuchbiggerthanlargestprimefactor}\end{equation}\label{dichotomywithpk}
\end{proposition}
\begin{proof} 

If $p_1=3$ and $k=3$, it is easy to check the short finite list that all of these satisfy at least one of the two bounds. We thus may assume that $(p_1, k) \neq (3,3)$.

As before, if $n$ is non-deficient we must have
$$H(n) = \prod_{i=1}^k \frac{p_i}{p_i-1} >2, $$
and thus 
$$\frac{p_1p_2p_3 \cdots p_k}{(p_1 -1)(p_2-1) \cdots (p_k-1)} -2 >0 $$
which means that \begin{equation}p_1p_2 \cdots p_k  = 2(p_1-1)(p_2-1)\cdots (p_k-1) +x \label{xintroduction} \end{equation} for some $x$. Note that $x$ is itself a positive integer. 
Now,  Equation \eqref{xintroduction} modulo $p_k$ becomes
$$0 \equiv 2(p_1-1)(p_2-1)\cdots(p_{k-1}-1)(-1) + x \pmod{p_k}.$$
So $p_k \mid -2(p_1-1)(p_2-1)\cdots(p_{k-1}-1) + x$.
If $-2(p_1-1)(p_2-1)\cdots(p_{k-1}-1) + x$ is negative, then
$p_k \mid 2(p_1-1)(p_2-1)\cdots(p_{k-1}-1) -x. $
Thus,
\begin{equation}p_k \leq 2(p_1-1)(p_2-1)\cdots(p_{k-1}-1) -x < 2(p_1-1)(p_2-1)\cdots(p_{k-1}-1) \label{initialpkinequalitytogetRform} \end{equation}
We have then 
\begin{equation} 2(p_1-1)(p_2-1)\cdots(p_{k-1}-1)= 2R\frac{p_1-1}{p_1}\frac{p_2-1}{p_2} \cdots \frac{p_{k-1}-1}{p_{k-1}}\frac{1}{p_k} \label{intermediateboundtogetRform}
\end{equation}

Since $n$ is non-deficient we have
$$\frac{p_1-1}{p_1}\frac{p_2-1}{p_2} \cdots \frac{p_{k-1}-1}{p_{k-1}}\frac{p_{k-1}}{p_{k}} <2, $$

and so 
\begin{equation} \frac{p_1-1}{p_1}\frac{p_2-1}{p_2} \cdots \frac{p_{k-1}-1}{p_{k-1}} < \frac{1}{2}\frac{p_k}{p_k-1} \label{withoutpkproductclosetoonehalf}.
\end{equation}

We then combine inequalities \eqref{initialpkinequalitytogetRform}, \eqref{intermediateboundtogetRform} and \eqref{withoutpkproductclosetoonehalf} to get that

$$p_k \leq \frac{2R}{p_k}\frac{1}{2}\frac{p_k}{p_k-1}$$ which implies that
\begin{equation}p_k (p_k-1) < R.\label{pktimespkminus1lessthanR}\end{equation}

Since $R$ is a multiple of $p_k$, we must have $$p_k(p_k-1) \leq R-p_k$$ which implies that $p_k < \sqrt{R}$.

Now, assume that $$p_k \geq  \sqrt{R}.$$ Thus, by our above reasoning we must that $-2(p_1-1)(p_2-1)\cdots(p_{k-1}-1) + x$  is non-negative. Then we have that $x \geq 2(p_1-1)(p_2-1)\cdots(p_{k-1}-1)$.

Thus, we have that
\begin{equation}p_1p_2 \cdots p_k  \geq  2(p_1-1)(p_2-1)\cdots (p_k-1) + 2(p_1-1)(p_2-1)\cdots(p_{k-1}-1) \label{Wegettousex},\end{equation}
and so
$$\prod_{i=1}^k \frac{p_i}{p_i-1} \geq 2 +\frac{2(p_1-1)(p_2-1)\cdots(p_{k-1}-1)}{(p_1-1)(p_2-1)\cdots(p_k-1)} =  2 + \frac{2}{p_k-1}.   $$

The rest of the proof then continues from invoking Lemma \ref{GeneralizedPuchtalemma} and proceeding as before. 
\end{proof}

One can use the same technique as in the proof of Proposition \ref{dichotomywithpk} to get a similar bound for even primitive non-deficient numbers. However, the result is weaker than the straightforward bound one gets just from analyzing how the power of 2 must behave per Lemma \ref{Classificationofprimitivenondeficientwhenkis2}.

Note that 945 satisfies Inequality \eqref{pklessthansquarerootofR} but not Inequality \eqref{mincomponentnotmuchbiggerthanlargestprimefactor}. In contrast, many primitive non-deficient numbers satisfy both inequalities. For example, 4095 satisfies both inequalities. Empirically it seems that most common behavior for odd non-deficient numbers is to satisfy both inequalities. An example of an odd primitive non-deficient number which does not satisfy Inequality  \eqref{pklessthansquarerootofR} is $$n = (3^2)(5^4)(11^3)(137^2)(31873),$$ for which we have $\frac{R^{\frac{1}{2}}}{p_k} \approx 0.84$. This raises the following question.

\begin{question} For all $\epsilon>0$, does there exist an odd, primitive non-deficient number $n=p_1^{a_1}p_2^{a_2} \cdots p_k^{a_k}$ with prime factors $p_1 < p_2 < \cdots < p_k$ such that
$$(\rad(n))^{\frac{1}{2}} < \epsilon p_k?$$  \label{questionaboutradicalandpk}
\end{question}

We suspect that the answer to Question \ref{questionaboutradicalandpk} is yes, but we are uncertain. At this point, we do not even see any way of proving that there are infinitely many odd primitive non-deficient numbers which do not satisfy  Inequality \eqref{pklessthansquarerootofR}.

It seems that most odd primitive non-deficient numbers satisfy both Inequality \eqref{pklessthansquarerootofR} and Inequality \eqref{mincomponentnotmuchbiggerthanlargestprimefactor}. It is not too hard to show that there are infinitely many which satisfy both. In particular, there are infinitely many square free odd primitive non-deficient numbers made from taking a prime and then multiplying by each successive prime until the number is not deficient. All numbers generated this way satisfy both inequalities.  Let $D$ be the set of odd primitive non-deficient numbers which satisfy both Inequality \eqref{pklessthansquarerootofR} and Inequality \eqref{mincomponentnotmuchbiggerthanlargestprimefactor}. Then we suspect the following. \begin{conjecture}
$$\lim_{x \rightarrow \infty} \frac{D(x)}{C_O(x)}=1.$$
\end{conjecture}

\section{Acknowledgments}The results in this paper were substantially tightened due to suggestions by Jan-Christoph Schlage-Puchta. The referee also made substantial improvements which greatly improved the presentation and helped fix gaps in the proofs of various lemmas as well as in Theorem \ref{maintheoremboundwithkR}, fixed the statement and proof of Proposition \ref{dichotomywithpk}, and pointed out McNew and Setty's recent work.

\bigskip
\hrule
\bigskip

\noindent 2020 {\it Mathematics Subject Classification}:
Primary 11A25; Secondary 11N64.

\noindent \emph{Keywords: } primitive non-deficient number, odd perfect number, abundant number

\bigskip
\hrule
\bigskip

\noindent (Concerned with sequences
\seqnum{A000396}, \seqnum{A005100}, and \seqnum{A006039}.)

\bigskip
\hrule
\bigskip

\vspace*{+.1in}
\noindent
Received 
revised version received  
Published in {\it Journal of Integer Sequences},

\bigskip
\hrule
\bigskip

\noindent
Return to \href{https://cs.uwaterloo.ca/journals/JIS/}{Journal of Integer Sequences home page}.
\vskip .1in


\begin{thebibliography}{99}
\bibitem{AcquaahKonyagin}  P. Acquaah and S. Konyagin, On prime factors of odd perfect numbers, \textit{Int. J. Number Theory} \textbf{8} (2012),  1537--1540.
\bibitem{AslaskenKirfel} H. Aslaksen and C. Kirfel,
Linear support for the prime number sequence and the first and second Hardy-Littlewood conjectures. 
\textit{Math. Scand.}\textbf{ 130} (2024),  198--220.
\bibitem{Avidon} M. Avidon,  On the distribution of primitive abundant numbers, {\it Acta Arith.}  \textbf{77} (1996), 195--205.
\bibitem{BVZ} S. Bibby, P.  Vyncke, and J. Zelinsky, 
On the third largest prime divisor of an odd perfect number, \textit{Integers} \textbf{21} (2021), \#A115.  
\bibitem{Billingsley} P. Billingsley,  On the distribution of large prime divisors. {\it Period. Math. Hungar.} \textbf{2} (1972), 283--289.
\bibitem{BKKL} K. Broughan, J. De Koninck, I. K\'{a}tai, and F. Luca, On integers for which the sum of divisors is the square of the squarefree core, \textit{J. Integer Seq.} 15 (2012), Article 12.7.5.
\bibitem{Deleglise} M. Del\'eglise,  Bounds for the density of abundant integers, {\it Experiment. Math.} \textbf{7}  (1998),  137--143.
\bibitem{Dickson} L.  Dickson,  Finiteness of the odd perfect and primitive abundant numbers with $n$ distinct prime factors, {\it Amer. J. Math.} \textbf{35} (1913), 413--422.
\bibitem{Ellia} P. Ellia, A remark on the radical of odd perfect numbers, 
\textit{Fibonacci Quart.} \textbf{50} (2012), 231--234.
\bibitem{OchemRaoradica} P. Ochem and M. Rao, 
Another remark on the radical of an odd perfect number, \textit{Fibonacci Quart.} \textbf{52} (2014), 215--217.
\bibitem{Grun}  O. Gr{\"u}n,  {\"U}ber ungerade vollkommene Zahlen, {\it Math. Z.} {\bf 55} (1952), 353--354.
\bibitem{Heath-Brown} D. Heath-Brown,  Odd perfect numbers,
\textit{Math. Proc. Cambridge Philos. Soc.} \textbf{115} (1994),  191--196.
\bibitem{Klurman} O. Klurman, 
Radical of perfect numbers and perfect numbers among polynomial values.
\textit{Int. J. Number Theory }\textbf{12} (2016), 585--591.
\bibitem{Kobayashi} M. Kobayashi,  A new series for the density of abundant numbers.  {\it Int. J. Number Theory} \textbf{10}  (2014),  73--84.
\bibitem{Luca} F. Luca, On numbers n for which the prime factors $\sigma(n)$ are among the prime factors of $n$, \textit{Result. Math.} \textbf{45} (2004), 79--87.
\bibitem{LucaPomerance} F. Luca and C. Pomerance, On the radical of a perfect number, \textit{New York J. Math.} \textbf{16} (2010), 23--30.
\bibitem{McNewSetty} N. McNew and J. Setty, On the densities of covering numbers and abundant numbers, arxiv preprint arXiv:2507.23041 [math.NT], 2025, available at {\url{https://arxiv.org/abs/2507.23041}}.
\bibitem{Nielsen} P. Nielsen, An upper bound for odd perfect numbers, \textit{Integers} \textbf{3} (2003), \#A14.
\bibitem{Servais} C. Servais, Sur les nombres parfaits, {\it Mathesis} {\bf 8} (1888), 92--93.
\bibitem{Stone}  A. Stone, Improved upper bounds for odd perfect numbers. I., {\it{Integers}} \textbf{24}, (2024), \#A114.
\bibitem{Yamada2} T. Yamada, On a problem of De Koninck, 
\textit{Mosc. J. Comb. Number Theory,} {\bf 10} (2021), 249--260. 
\bibitem{Zelinskybig}  J. Zelinsky, On the total number of prime factors of an odd perfect number, {\it Integers} {\bf 21} (2021), \#A76.

\end{thebibliography}
\end{document}